\documentclass{amsart}
\date{\today}
\usepackage{amssymb}
\usepackage{latexsym}
\newcommand{\Z}{{\mathbb Z}}
\newcommand{\R}{{\mathbb R}}

\newcommand{\beq}{\begin{eqnarray}}
\newcommand{\eeq}{\end{eqnarray}}

\newtheorem{theorem}{Theorem}

\newtheorem{lemma}{Lemma}[section]

\begin{document}

\title[Unbounded Jacobi Matrices at Critical Coupling]{The Case of Critical
Coupling in a Class of Unbounded Jacobi Matrices Exhibiting a First-Order Phase
Transition}

\author[D.\ Damanik, S.\ Naboko]{David Damanik$\, ^1$ and
Serguei Naboko$\, ^2$}

\thanks{D.\ D.\ was supported in part by NSF grants DMS--0227289 and DMS--0500910.\\
\hspace*{9.5pt} S.\ N.\ was supported in part by RFBR 06-01-00249 and would also like to
express his gratitude to the Mathematics Department of Caltech where most of this work
was done.}

\maketitle
\vspace{0.3cm}
\noindent $^1$ Department of Mathematics
253--37, California Institute of Technology, Pasadena,
CA 91125, USA, E-mail: \mbox{damanik@its.caltech.edu}\\[2mm]
$^2$ Department of Mathematical Physics, Institute of Physics,
St.~Petersburg University, Ulianovskaia~1, 198904~St.~Petergoff,
St.~Petersburg, Russia,
E-mail: \mbox{naboko@snoopy.phys.spbu.ru}\\[3mm]

\begin{abstract}
We consider a class of Jacobi matrices with unbounded
coefficients. This class is known to exhibit a first-order phase
transition in the sense that, as a parameter is varied, one has
purely discrete spectrum below the transition point and purely
absolutely continuous spectrum above the transition point. We
determine the spectral type and solution asymptotics at the
transition point.
\end{abstract}

%
%
%
%

\section{Introduction}

In this paper we analyze spectral properties of Jacobi matrices,
\begin{equation}\label{jacmatint}
J = \begin{pmatrix}
b_1 & a_1 & 0 & 0 & \cdots \\
a_1 & b_2 & a_2 & 0 & \cdots \\
0 & a_2 & b_3 & a_3 & \cdots \\
\vdots & \vdots & \vdots & \vdots & \ddots
\end{pmatrix},
\end{equation}
acting in $\ell^2 (\Z_+)$, and asymptotic properties of solutions
to the associated difference equation
\begin{equation}\label{eveint}
a_n u_{n+1} + b_n u_n + a_{n-1} u_{n-1} = E u_n.
\end{equation}

Recently, there has been interest in the case of unbounded
coefficients $a_n,b_n$; see, for example,
\cite{d,dp1,dp2,jm,jn1,jn2,jn3,jns,m}.

Motivated in particular by \cite{e2,lln1,lln2,mr}, Janas and
Naboko \cite{jn2} studied a large class of Jacobi matrices with
unbounded and periodically modulated entries. To be specific, they
considered the case $a_n = c_n \mu_n$, $b_n = d_n r_n$, where
$\{c_n\}$ is strictly positive and $N$-periodic and $\{d_n\}$ is
$M$-periodic. The sequences $\{\mu_n\}$, $\{r_n\}$ are unbounded
and satisfy a number of conditions. The main example one should
have in mind is $\mu_n = r_n = n^\alpha$, where $0 < \alpha \le
1$. Let us consider this special case, that is,
\begin{equation}\label{ex}
a_n = c_n n^\alpha, \; \; \; b_n = d_n n^\alpha.
\end{equation}
Since the Carleman condition holds, that is,
$$
\sum_{n=1}^\infty a_n^{-1} = \infty,
$$
$J$ defines a self-adjoint operator in $\ell^2 (\Z_+)$. The
spectral type of $J$ is closely related to the location of zero
relative to the spectrum of the associated periodic Jacobi matrix
$J_{{\rm per}}$ which is given by
$$
J_{{\rm per}} = \begin{pmatrix}
d_1 & c_1 & 0 & 0 & \cdots \\
c_1 & d_2 & c_2 & 0 & \cdots \\
0 & c_2 & d_3 & c_3 & \cdots \\
\vdots & \vdots & \vdots & \vdots & \ddots
\end{pmatrix}.
$$
Clearly, $J_{{\rm per}}$ is $K$-periodic, where $K$ is the least
common multiple of $M$ and $N$. Its characteristic polynomial,
$d(E)$, is given by
$$
d(E) = \text{Tr} \left[ \prod_{n=1}^K \left( \begin{array}{cc} 0 &
1 \\ - c_{n-1} c_n^{-1} & (E - d_n) c_n^{-1} \end{array} \right)
\right].
$$
It is well known that the spectrum of $J_{{\rm per}}$ is given by
$$
\sigma(J_{{\rm per}}) = \{ E \in \R : |d(E)| \le 2 \}.
$$
This set is the union of $K$ bands (non-degenerate closed
intervals) whose interiors are mutually disjoint. The following
was shown in \cite{jn2}:

\begin{theorem}[Janas-Naboko]\label{jannabthm}
{\rm (a)} If $|d(0)| < 2$, then the spectrum of $J$ is purely
absolutely continuous and $\sigma(J) = \R$.
\\[1mm]
{\rm (b)} If $|d(0)| > 2$, then the spectrum of $J$ is purely
discrete.
\end{theorem}

In other words, if zero is not one of the band edges of
$\sigma(J_{{\rm per}})$, one has a complete understanding of the
spectral type of $J$. Moreover, Janas and Naboko also described
the asymptotic behavior of the solutions of \eqref{eveint};
compare \cite[Theorems~3.1 and 4.2]{jn2}.

The question of what happens at transition points, corresponding
to $d(0) = \pm 2$, was left open in \cite{jn2}. It was pointed out
that new methods and ideas would be necessary to understand these
critical cases. It is our goal here to study this scenario in a
simple special case. We shall see that even in this simple
situation, the analysis is already quite involved. Moreover, in
this way the main new ideas are more transparent.

We will study the case where the $a_n$'s and $b_n$'s are given by
\eqref{ex} and both periods, $M$ and $N$, are small. Specifically,
let us consider the case
$$
M=2, \;\;\; N=1, \;\;\; c_n \equiv 1, \;\;\; d_{2n-1} \equiv b, \;\;\; d_{2n} \equiv
\tilde{b}.
$$
We find
$$
d(0) = -2 + b\tilde{b}.
$$
Fix $b > 0$. Then, by Theorem~\ref{jannabthm}, $J$ has purely
absolutely continuous spectrum for (not too large) $\tilde{b} > 0$
and purely discrete spectrum for $\tilde{b} < 0$. Similarly, for
$b < 0$, $J$ has purely absolutely continuous spectrum for
$\tilde{b} < 0$ (again, the absolutely value should not be too
large) and purely discrete spectrum for $\tilde{b}
> 0$. Thus, if we fix some non-zero value for $b$, the case
$\tilde{b} = 0$ is the critical case for which
Theorem~\ref{jannabthm} does not apply.

It is our goal to study this particular case, that is, we want to
determine the spectral type and solution asymptotics for
$$
a_n = n^\alpha, \;\; b_n = \left\{ \begin{array}{cl} bn^\alpha &
\text{ if $n$ is odd,} \\ 0 & \text{ if $n$ is even,} \end{array}
\right.
$$
where $b \not= 0$ and $0 < \alpha \le 1$. Our main results in the
case $b > 0$ (when $b < 0$, one has to reflect the energy about
zero) are as follows:
\begin{itemize}

\item If $2/3 < \alpha \le 1$, the spectrum of $J$ is purely
absolutely continuous on $(-\infty,0)$. Moreover, explicit
solution asymptotics are given. These results are stated in more
detail and proven in Section~\ref{three}.

\item If $0 < \alpha \le 1$, zero is not an eigenvalue of $J$. See
Section~\ref{four}.

\item If $0 < \alpha \le 1$, the spectrum of $J$ is purely
discrete in $(0,\infty)$. The eigenvalues are simple and the
$n$-th eigenvalue, $E_n$, obeys the bounds
$$
C_1(b) n^\alpha \le E_n \le C_2(b) n^\alpha.
$$
We also provide explicit expressions for the (positive and finite)
constants $C_1(b), C_2(b)$. These results and their proofs can be
found in Section~\ref{five}.

\end{itemize}

This determines the spectral type completely when $2/3 < \alpha \le 1$. In this case, the
spectrum looks essentially like

\medskip

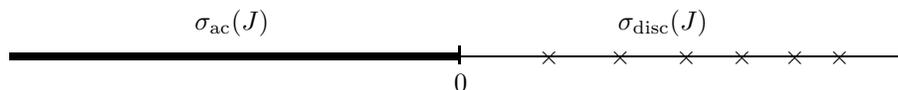
\begin{picture}(340,40)

\put(0,13){\line(1,0){340}}

\put(170,10){\line(0,1){7}}

\put(168,0){0}

\put(200,10){$\times$}

\put(227,10){$\times$}

\put(252,10){$\times$}

\put(273,10){$\times$}

\put(293,10){$\times$}

\put(310,10){$\times$}

\put(70,23){$\sigma_{{\rm ac}}(J)$}

\put(230,23){$\sigma_{{\rm disc}}(J)$}

\linethickness{1mm}

\put(0,13){\line(1,0){170}}

\end{picture}

\medskip

The condition $2/3 < \alpha \le 1$ is naturally associated with our method of proof
(compare, e.g., \eqref{alphacondition}). However, we expect the picture above also when
$0 < \alpha \le 2/3$. Thus, we leave the question of proving purely absolutely continuous
spectrum on the negative energy axis for these values of $\alpha$ as an open problem.
Other open problems suggested by our work will be discussed in Section~\ref{six}.

\medskip

The spectral analysis in the energy region $(-\infty,0]$ is based on an analysis of the
solutions to the difference equation \eqref{eveint}. We obtain asymptotic expressions for
all solutions corresponding to energies $E \in (-\infty,0]$; see
Theorems~\ref{acspectrum} and \ref{zeroenergyasymptotics}. In particular, this determines
the asymptotic behavior of the orthogonal polynomials associated with the spectral
measure of the pair $(J,\delta_1)$ since, by standard theory, they solve \eqref{eveint}
with a Dirichlet boundary condition, $u_0 = 0$.

%
%
%
%

\section{Preliminaries}\label{two}

We will study the Jacobi matrix
\begin{equation}\label{jacmat}
J = \begin{pmatrix}
b_1 & a_1 & 0 & 0 & \cdots \\
a_1 & b_2 & a_2 & 0 & \cdots \\
0 & a_2 & b_3 & a_3 & \cdots \\
\vdots & \vdots & \vdots & \vdots & \ddots
\end{pmatrix},
\end{equation}
acting in $\ell^2 (\Z_+)$, where the parameters $a_n, b_n$ are
given by
\begin{equation}\label{parameters}
a_n = n^\alpha, \;\; b_n = \left\{ \begin{array}{cl} bn^\alpha &
\text{ if $n$ is odd,} \\ 0 & \text{ if $n$ is even.} \end{array}
\right.
\end{equation}
Here, $b \not=0$. Note that if $J_{\{a_n\},\{b_n\}}$ denotes the
Jacobi matrix corresponding to the sequences $\{a_n\}$ and
$\{b_n\}$, and $U$ denotes the unitary transformation of $\ell^2$,
given by $(U \psi)_n = (-1)^n \psi_n$, then
$$
U J_{\{a_n\},\{b_n\}} U = - J_{\{a_n\},\{-b_n\}}.
$$
We will therefore restrict our attention in what follows to the
case $b > 0$. Consider solutions of the difference equation
\begin{equation}\label{eve}
a_n u_{n+1} + b_n u_n + a_{n-1} u_{n-1} = E u_n.
\end{equation}
Defining
$$
U_n = \left( \begin{array}{c} u_{n-1} \\ u_n \end{array} \right),
$$
the recursion \eqref{eve} is equivalent to
$$
U_{n+1} = T_n U_n,
$$
where
$$
T_n = \left( \begin{array}{cc} 0 & 1 \\ -\frac{a_{n-1}}{a_n} &
\frac{E - b_n}{a_n} \end{array} \right).
$$
Let
$$
B_n = T_{2n} T_{2n-1},
$$
so that
\begin{equation}\label{eve2}
U_{2n+1} = (B_n \times \cdots \times B_1) U_1.
\end{equation}
We have
\begin{align*}
B_n & = \left( \begin{array}{cc} 0 & 1 \\ -(1 -
\frac{1}{2n})^\alpha & \frac{E}{(2n)^\alpha}
\end{array} \right) \left( \begin{array}{cc} 0 & 1 \\ -(1 -
\frac{1}{2n-1})^\alpha & \frac{E - b(2n-1)^\alpha}{(2n-1)^\alpha} \end{array} \right) \\
& = \left( \begin{array}{cc} -1 & -b \\ 0 & -1
\end{array} \right) + \frac{1}{(2n)^\alpha} \left( \begin{array}{rc} 0 & E \\ -E &
-bE \end{array} \right) + \frac{1}{2n} \left( \begin{array}{cc}
\alpha & 0 \\ 0 & \alpha \end{array} \right) + O(n^{-2\alpha}).
\end{align*}
We see from \eqref{eve2} that we should study the product of the
form $B_n \times \cdots \times B_1$. However, we shall study a
slightly different sequence of products whose usefulness will
become clearer below. To define this auxiliary problem, we have to
introduce a few matrices. Let
$$
C_n = \left( \begin{array}{rr} 1 & -b \\ 1 & 0
\end{array} \right) + \frac{1}{(2n)^\alpha} \left( \begin{array}{cc} bE + E/(2b) & 0 \\
E/(2b) & -E/2 \end{array} \right) + \frac{1}{2n} \left(
\begin{array}{rc} 0 & 0 \\ -\alpha & \alpha b \end{array} \right)
$$
and
$$
\tilde{B}_n = \left( \begin{array}{rc} 0 & 1 \\ -1 & 2
\end{array} \right) + \frac{1}{(2n)^\alpha} \left( \begin{array}{cc} 0 & 0 \\ 0 &
bE \end{array} \right) + \frac{1}{n} \left( \begin{array}{cr} 0 &
0
\\ \alpha & -\alpha \end{array} \right).
$$

\begin{lemma}\label{L:conjugacy}
We have
$$
C_n B_n C_n^{-1} = - \tilde{B}_n + O(n^{-2\alpha}).
$$
\end{lemma}

\begin{proof}
This is tedious but straightforward.
\end{proof}

%
%
%
%

\section{The Absolutely Continuous Spectrum and Solution Asymptotics at Negative
Energies}\label{three}

Motivated by \eqref{eve2} and Lemma~\ref{L:conjugacy}, we will
study asymptotics for the auxiliary problem
\begin{equation}\label{eve3}
V_{n+1} = (\tilde{B}_n \times \cdots \times \tilde{B}_1) V_1.
\end{equation}
This problem is a more general version of the one studied in
\cite[Section~3]{jn1}. There, products of matrices of the form
$$
\hat{B}_n = \left( \begin{array}{rc} 0 & 1 \\ -1 & 2
\end{array} \right) + \frac{1}{n} \left( \begin{array}{cc} 0 &
0 \\ 1 & \lambda \end{array} \right)
$$
were studied. We will employ a similar strategy to find solution
asymtotics for the more general problem at hand. This explains why
we introduced the matrices $\tilde{B}_n$ and derived
Lemma~\ref{L:conjugacy} in the previous section.

\begin{theorem}\label{negensolas}
Suppose that $2/3 < \alpha \le 1$ and $b > 0$. Then, for every $E
< 0$, \eqref{eve3} has two linearly independent solutions
\begin{equation}\label{thm1-1}
V_n^\pm = \left( \begin{array}{c} v_{n-1}^\pm \\ v_n^\pm
\end{array} \right)
\end{equation}
with asymptotics given by
\begin{equation}\label{thm1-2}
v_n^\pm = n^{-\alpha/4} \exp \left( \pm i
\frac{\sqrt{-bE}}{2^{\alpha/2}} \cdot \frac{n^{1 -
\frac{\alpha}{2}}}{1 - \frac{\alpha}{2}} \right) (1 + o(1)) \text{
as } n \to \infty.
\end{equation}
\end{theorem}

\begin{proof}
We make the ansatz
$$
z_n = n^\gamma \exp \left( A n^\delta \right)
$$
and define the matrix
$$
S_n = \left( \begin{array}{cc} \overline{z_{n-1}} & z_{n-1} \\
\overline{z_n} & z_n \end{array} \right).
$$
Our goal is to choose $\gamma, A, \delta$ such that
\begin{equation}\label{zngoal}
S_{n+1}^{-1} \tilde{B}_n S_n = I + R_n, \; \{ \| R_n \| \} \in
\ell^1.
\end{equation}
Consequently, an arbitrary non-trivial solution of \eqref{eve3}
has the form $V_n = S_n W_n$, where $W_n$ is a sequence of vectors
which tends to a non-zero vector $W$. Since \eqref{zngoal} will be
shown to hold if we let
\begin{equation}\label{zndef}
\gamma = - \frac{\alpha}{4}, \; A = \frac{\sqrt{bE}}{2^{\alpha/2}
(1 - \frac{\alpha}{2})}, \; \delta = 1 - \frac{\alpha}{2},
\end{equation}
the assertion of the theorem then follows immediately.

Let us consider the matrix $S_{n+1}^{-1} \tilde{B}_n S_n$. It is
readily checked that
$$
S_{n+1}^{-1} \left( \begin{array}{rc} 0 & 1 \\
-1 & 2 \end{array} \right) S_n = (\det S_{n+1})^{-1} \left(
\begin{array}{rr}
x_n & y_n \\
-\overline{y_n} & -\overline{x_n} \end{array} \right),
$$
where
\begin{align*}
x_n & = z_n \overline{z_{n-1}} + \overline{z_n} z_{n+1} -
2|z_n|^2,\\
y_n & = z_n z_{n-1} + z_n z_{n+1} - 2 z_n^2,
\end{align*}
and
$$
S_{n+1}^{-1} \left( \begin{array}{rr} 0 & 0 \\
a & b \end{array} \right) S_n = (\det S_{n+1})^{-1} \left(
\begin{array}{rr}
s_n & t_n \\
-\overline{t_n} & -\overline{s_n} \end{array} \right),
$$
where
\begin{align*}
s_n & = - a z_n \overline{z_{n-1}} -
b |z_n|^2,\\
t_n & = - a z_n z_{n-1} - b z_n^2.
\end{align*}
Putting this together, we obtain
$$
S_{n+1}^{-1} \tilde{B}_n S_n = (\det S_{n+1})^{-1} \left(
\begin{array}{rr}
a_n & b_n \\
-\overline{b_n} & -\overline{a_n} \end{array} \right),
$$
where
\begin{align*}
a_n & = z_n \overline{z_{n-1}} + \overline{z_n} z_{n+1} - 2|z_n|^2
+ \tfrac{1}{(2n)^\alpha} \left[ - bE |z_n|^2 \right] +
\tfrac{1}{n}
\left[ - \alpha z_n \overline{z_{n-1}} + \alpha |z_n|^2 \right],\\
b_n & = z_n z_{n-1} + z_n z_{n+1} - 2 z_n^2 +
\tfrac{1}{(2n)^\alpha} \left[ - bE z_n^2 \right] + \tfrac{1}{n}
\left[ - \alpha z_n z_{n-1} + \alpha z_n^2 \right] .
\end{align*}

We consider first the off-diagonal elements of $S_{n+1}^{-1}
\tilde{B}_n S_n$:
\begin{align*}
(\det S_{n+1})^{-1} b_n & = \frac{z_n z_{n-1} + z_n z_{n+1} - 2
z_n^2 - \tfrac{bE}{(2n)^\alpha} z_n^2 + \tfrac{\alpha}{n} \left[
z_n^2 - z_n z_{n-1} \right]}{z_{n+1} \overline{z_n} - \overline{z_{n+1}} z_n} \\
& = \frac{z_n^{-1} z_{n-1} + z_n^{-1} z_{n+1} - 2 -
\tfrac{bE}{(2n)^\alpha} + \tfrac{\alpha}{n} \left[ 1 - z_n^{-1}
z_{n-1} \right]}{z_{n+1} \overline{z_n} z_n^{-2} -
\overline{z_{n+1}} z_n^{-1}}.
\end{align*}

For the denominator in the last expression, we find
$$
z_{n+1} \overline{z_n} z_n^{-2} - \overline{z_{n+1}} z_n^{-1} = 2A
\delta n^{\delta - 1} (1 + o(1)),
$$
provided we choose $\delta$ as in \eqref{zndef} (which implies
$1/2 \le \delta \le 2/3$).

For the numerator, we find after lengthy but straightforward
calculations,
$$
z_n^{-1} z_{n-1} + z_n^{-1} z_{n+1} - 2 - \tfrac{bE}{(2n)^\alpha}
+ \tfrac{\alpha}{n} - \tfrac{\alpha}{n} z_n^{-1} z_{n-1} = X
n^{\delta - 2} + Y n^{2\delta - 2} + O(n^{4 \delta - 4}),
$$
where
$$
X = A \delta (\delta - 1) + 2 \gamma A \delta + \alpha A \delta,
\; \; Y = (A \delta)^2 - \frac{bE}{2^\alpha}.
$$
If $\gamma, A, \delta$ are chosen as in \eqref{zndef}, then $X = Y
= 0$ and hence
\begin{equation}\label{alphacondition}
(\det S_{n+1})^{-1} b_n = O \left( n^{3\delta - 3} \right) = O
\left( n^{- \frac{3}{2} \alpha} \right) ,
\end{equation}
which is summable by our assumption on $\alpha$.

Next, we consider the diagonal elements of $S_{n+1}^{-1}
\tilde{B}_n S_n$:
\begin{align*}
| a_n - \det S_{n+1}| & = \left| z_n \overline{z_{n-1}} +
\overline{z_{n+1}} z_n - 2|z_n|^2 - \tfrac{bE
|z_n|^2}{(2n)^\alpha} + \tfrac{\alpha}{n} \left[ - z_n
\overline{z_{n-1}} + |z_n|^2 \right] \right| \\
& = \left| \overline{z_n} z_{n-1} + z_{n+1} \overline{z_n} -
2|z_n|^2 - \tfrac{bE |z_n|^2}{(2n)^\alpha} + \tfrac{\alpha}{n}
\left[ - \overline{z_n} z_{n-1} + |z_n|^2 \right] \right| \\
& = \left| z_n z_{n-1} + z_n z_{n+1} - 2 z_n^2 - \tfrac{bE
z_n^2}{(2n)^\alpha} + \tfrac{\alpha}{n} \left[ - z_n z_{n-1} +
z_n^2 \right] \right| \\
& = | b_n |,
\end{align*}
and hence $\left| (\det S_{n+1})^{-1} a_n - 1 \right| = \left|
(\det S_{n+1})^{-1} b_n \right|$. This shows that \eqref{zngoal}
holds, concluding the proof.
\end{proof}

The special case $\alpha = 1$ deserves an additional remark.
Theorem~\ref{negensolas} above gives the asymptotics
\begin{equation}\label{alphaeqone}
n^{-1/4} \exp \left( \pm i \sqrt{-2 b E n} \right) (1 + o(1))
\text{ as } n \to \infty
\end{equation}
for a pair of linearly independent solutions of the auxiliary
problem \eqref{eve3}. This can also be shown by an application of
the Birkhoff-Adams Theorem; see \cite[Theorem~8.36]{e}. The latter
result concerns second-order difference equations of the form
\begin{equation}\label{gendiffeq}
x(n+2) + p_1(n) x(n+1) + p_2(n) x(n) = 0,
\end{equation}
where $p_1(n)$ and $p_2(n)$ have asymptotic expansions
$$
p_1(n) = \sum_{j=0}^\infty \frac{c_j}{n^j} , \;\;\; p_2(n) =
\sum_{j=0}^\infty \frac{d_j}{n^j}
$$
with $d_0 \not=0$. In our situation, we have
$$
c_0 = -2, \; \; \; c_1 = 1 - \frac{b E}{2}, \; \; \; d_0 = 1, \;
\; \; d_1 = - 1.
$$
We have to compute the roots $\lambda_1, \lambda_2$ of the
characteristic equation $\lambda^2 + c_0 \lambda + d_0 = 0$ and
find $\lambda_1 = \lambda_2 = 1$. Moreover, if $E \not= 0$, we
have $2d_1 \not= c_0 c_1$. Thus, part (b) of
\cite[Theorem~8.36]{e} tells us that there are two linearly
independent solutions $x_1(n)$, $x_2(n)$ of \eqref{gendiffeq}
whose asymptotics are given by \eqref{alphaeqone}. (Note, however,
that there is a misprint in formula (8.6.7) of \cite{e}:
$n^{\gamma_i}$ should be replaced by $n^\alpha$.)

\begin{theorem}\label{acspectrum}
Suppose that $2/3 < \alpha \le 1$ and $b > 0$. Then, for every $E
< 0$, \eqref{eve2} has two linearly independent solutions
$$
U_{2n+1}^\pm = \left( \begin{array}{c} u_{2n}^\pm \\ u_{2n+1}^\pm
\end{array} \right)
$$
with asymptotics given by
$$
U_{2n+1}^\pm = (-1)^n T V_n^\pm (1 + o(1)),
$$
where
$$
T = \frac{1}{b} \left( \begin{array}{rr} 0 & b \\ -1 & 1
\end{array} \right)
$$
and $V_n^\pm$ are given by \eqref{thm1-1} and \eqref{thm1-2}.
Moreover, the spectrum of $J$ is purely absolutely continuous on
$(-\infty,0)$.
\end{theorem}

\begin{proof}
By standard results of asymptotic analysis (centered around
Levinson's Theorem; compare \cite[Section~8.3]{e} and
\cite{Eastham}), Lemma~\ref{L:conjugacy} implies that the
solutions of \eqref{eve2} behave asymptotically like the solutions
of the system
$$
\tilde{U}_{2n+1} = \left[ (-C_n^{-1} \tilde{B}_n C_n) \times
\cdots \times (-C_1^{-1} \tilde{B}_1 C_1) \right] \tilde{U}_1.
$$
Since $C_m C_{m-1}^{-1} = I + O(m^{-2})$, we may as well consider
the system
$$
\hat{U}_{2n+1} = (-1)^n C_n^{-1} \left[ \tilde{B}_n  \times \cdots
\times \tilde{B}_1 \right] \hat{U}_1.
$$
Now the first assertion is an immediate consequence of
Theorem~\ref{negensolas}.

To prove that the spectrum of $J$ is purely absolutely continuous
on $(-\infty,0)$, it suffices to show that for every $E < 0$,
there does not exist a subordinate solution of \eqref{eve} in the
sense of Gilbert-Pearson (cf.~\cite{gp}; see also \cite{jl,kp}).
Thus, we claim that for every pair of non-trivial solutions
$u^{(1)}$, $u^{(2)}$ of \eqref{eve}, we have
\begin{equation}\label{nonsubord}
\limsup_{N \to \infty} \left[ \frac{\sum_{n=1}^N \left| u^{(1)}_n
\right|^2}{\sum_{n=1}^N \left| u^{(2)}_n \right|^2} \right] > 0.
\end{equation}
Since every solution $u$ of \eqref{eve} is given by a linear
combination of the two solutions whose asymptotic behavior is
given by $(-1)^n T V_n^\pm (1 + o(1))$, there are ($u$-dependent)
positive constants $C_1, C_2$ such that
$$
C_1 N^{1 - \frac{\alpha}{2}} \le \sum_{n=1}^N \left| u_n \right|^2
\le C_2 N^{1 - \frac{\alpha}{2}}.
$$
From this, \eqref{nonsubord} follows immediately.
\end{proof}

%
%
%
%

\section{Solution Asymptotics at Zero Energy}\label{four}

In this section we study the asymptotics of the matrix product
$B_n \times \cdots \times B_1$, and hence the asymptotics of
solutions, for $E = 0$. We show in particular that there are no
non-trivial solutions in $\ell^2$.

\begin{theorem}\label{zeroenergyasymptotics}
For zero energy, $E = 0$, the difference equation \eqref{eve} has
two linearly independent solutions $u_n^{(j)}$, $j = 1,2$, with
asymptotics given by
$$
\left( \begin{array}{l} u_{2n}^{(1)} \\ u_{2n+1}^{(1)} \end{array}
\right) = \left( \begin{array}{c} (-1)^n n^{-\frac{\alpha}{2}} (1
+ o(1))
\\ 0 \end{array} \right)
$$
and
$$
\left( \begin{array}{l} u_{2n}^{(2)} \\ u_{2n+1}^{(2)} \end{array}
\right) = \left( \begin{array}{c} (-1)^n \left[ b
n^{1-\frac{\alpha}{2}} + O(n^{- \frac{\alpha}{2}} \log n) \right]
(1 + o(1))
\\ (-1)^n n^{-\frac{\alpha}{2}} (1 + o(1))
\end{array} \right).
$$
In particular, zero is not an eigenvalue of $J$.
\end{theorem}

\begin{proof}
We have
\begin{align*}
B_n & = \left( \begin{array}{cc} 0 & 1 \\ -(1 -
\frac{1}{2n})^\alpha & 0
\end{array} \right) \left( \begin{array}{cc} 0 & 1 \\ -(1 -
\frac{1}{2n-1})^\alpha & - b \end{array} \right) \\
& = \left( \begin{array}{cc} -1 & -b \\ 0 & -1
\end{array} \right) +  \left( \begin{array}{cc}
\frac{\alpha}{2n} & 0 \\ 0 & \frac{\alpha}{2n} \end{array} \right)
+ O(n^{-2}) \left( \begin{array}{cc} r_n & 0 \\ 0 & s_n
\end{array} \right)
\end{align*}
with bounded sequences $\{r_n\}$, $\{s_n\}$.

That is,
$$
B_n = - M_b + \frac{\alpha}{2n} I + \left( \begin{array}{cc}
O(n^{-2}) & 0 \\ 0 & O(n^{-2})
\end{array} \right), \text{ where }
M_b = \left( \begin{array}{cc} 1 & b \\ 0 & 1
\end{array} \right).
$$
Therefore,
$$
B_n \times \cdots \times B_1 = (-1)^n M_b^n \left( M_b^{-n} (-B_n)
M_b^{n-1} \right) \times \cdots \times \left( M_b^{-1} (-B_1)
M_b^0 \right).
$$
Since
$$
M_b^{-j} (-B_j) M_b^{j-1} = I - \frac{\alpha}{2j} \left(
\begin{array}{cr} 1 & -b \\ 0 & 1
\end{array} \right) + \left( \begin{array}{cc}
O(j^{-2}) & O(j^{-1}) \\ 0 & O(j^{-2})
\end{array} \right),
$$
we obtain
$$
B_n \times \cdots \times B_1 = (-1)^n \left( \begin{array}{cc} 1 &
bn \\ 0 & 1 \end{array} \right) \left[ \prod_{j=1}^n \left( 1 -
\frac{\alpha}{2j} \right) \right] U_n,
$$
where
\begin{align*}
U_n & = \prod_{j=1}^n \left[ \left(
\begin{array}{cc} 1 + O(j^{-2}) & 0 \\ 0 & 1 + O(j^{-2})
\end{array} \right) + \left( \begin{array}{cc} 0 & O(j^{-1}) \\ 0 &
0 \end{array} \right) \right] \\
& = \left( \begin{array}{cc} C_1 + o(1) & O(\log n) \\ 0 & C_2 +
o(1)
\end{array} \right).
\end{align*}
Here, $C_1, C_2$ are suitable non-zero constants.

On the other hand,
$$
\prod_{j=1}^n \left( 1 - \frac{\alpha}{2j} \right) =
n^{-\frac{\alpha}{2}} ( C_3 + o(1) )
$$
with a suitable constant $C_3 \not= 0$.

Putting everything together, we find that
$$
B_n \times \cdots \times B_1 = (-1)^n n^{-\frac{\alpha}{2}} \left(
\begin{array}{cc} C_1 C_3 + o(1) & C_2 C_3 bn + O(\log n) \\ 0 & C_2 C_3 + o(1)
\end{array} \right).
$$
The assertion follows.
\end{proof}

%
%
%
%

\section{The Discrete Spectrum}\label{five}

In this section we study positive energies and prove that the
spectrum in $(0,\infty)$ is purely discrete and accumulates only
at $\infty$. Moreover, we provide upper and lower bounds for
eigenvalues.

Fix the parameters $b > 0$ and $0 < \alpha \le 1$. For $a > 0$, we
define
$$
\mathcal{H}^{(1)}_a = \text{span} \left\{ \delta_{2n - 1} : 1 \le
n \le \frac12 \left( \frac{2a}{b} \right)^{1/\alpha} \right\}, \;
\; \mathcal{H}^{(2)}_a = \left( \mathcal{H}^{(1)}_a \right)^\bot.
$$

\begin{lemma}\label{l:keyest}
For $a > 0$ and $\psi \in \mathcal{H}^{(2)}_a \cap D(J)$, we have
$ \left\| (J - aI) \psi \right\| \ge a \| \psi \|$.
\end{lemma}

\begin{proof}
Consider the even and odd subspaces of $\ell^2$,
$$
\ell^2_{{\rm e}} = \{ f \in \ell^2 : f_{2n-1} = 0 \; \forall n \in
\Z_+ \}, \; \; \ell^2_{{\rm o}} = \{ f \in \ell^2 : f_{2n} = 0 \;
\forall n \in \Z_+ \}
$$
and the respective orthogonal projections, $P_{{\rm e}}$ and
$P_{{\rm o}}$. We write a vector $\psi \in \ell^2$ as $\psi =
\psi_{{\rm e}} + \psi_{{\rm o}} = P_{{\rm e}} \psi + P_{{\rm o}}
\psi$.

Let $S$ be the shift, given by $(S \psi)_n = \psi_{n+1}$, and let
$A$ and $B$ be diagonal matrices with $A_{n,n} = a_n$ and $B_{n,n}
= b_n$, respectively. Then we can write the Jacobi matrix $J$ in
the form
$$
J = SA + AS^* + B.
$$
If we write
$$
J - aI = \left( \begin{array}{c|c} P_{{\rm o}} (J - aI) P_{{\rm
o}} & P_{{\rm o}} (J - aI) P_{{\rm e}} \\ \hline P_{{\rm e}} (J -
aI) P_{{\rm o}} & P_{{\rm e}} (J - aI) P_{{\rm e}} \end{array}
\right)
$$
and use
$$
B P_{{\rm e}} = 0, \; P_{{\rm o}} (SA + AS^*) P_{{\rm o}} = 0, \;
P_{{\rm e}} (SA + AS^*) P_{{\rm e}} = 0,
$$
we obtain
$$
J - aI = \left( \begin{array}{c|c} (B - aI) P_{{\rm o}} & P_{{\rm
o}} J_A P_{{\rm e}} \\ \hline P_{{\rm e}} J_A P_{{\rm o}} & - a
P_{{\rm e}} \end{array} \right).
$$
Here, $J_A$ denotes the matrix $SA + AS^*$. This yields
$$
(J - aI)^2 - a^2 I = \left( \begin{array}{c|c} (B^2 - 2aB) P_{{\rm
o}} + ( P_{{\rm e}} \ J_A P_{{\rm o}} )^* ( P_{{\rm e}} J_A
P_{{\rm o}} ) & (B - aI) P_{{\rm o}} J_A P_{{\rm e}}
- a P_{{\rm o}} J_A P_{{\rm e}} \\
\hline P_{{\rm e}} J_A P_{{\rm o}} (B - aI) P_{{\rm o}} - a
P_{{\rm e}} J_A P_{{\rm o}} & (P_{{\rm o}} J_A P_{{\rm e}})^*
(P_{{\rm o}} J_A P_{{\rm e}})
\end{array} \right)
$$

We want to show that
\begin{equation}\label{lemmagoal}
\langle [(J - aI)^2 - a^2 I] \psi, \psi \rangle \ge 0 \text{ for }
\psi \in \mathcal{H}^{(2)}_a.
\end{equation}

Observe first that
\begin{equation}\label{step1}
(B^2 - 2aB) P_{{\rm o}} \ge 0 \text{ on } \mathcal{H}^{(2)}_a,
\end{equation}
since
$$
(b (2n-1)^\alpha)^2 - 2a b (2n-1)^\alpha \ge 0 \Leftrightarrow
2n-1 \ge \left( \frac{2a}{b} \right)^{1/\alpha}.
$$

Now, let $\psi \in \mathcal{H}^{(2)}_a \cap D(J)$ be given and
write
$$
\psi_{{\rm e}} = P_{{\rm e}} \psi, \; \psi_{{\rm o}} = P_{{\rm o}}
\psi, \; \phi_{{\rm e}} = P_{{\rm e}} J_A \psi_{{\rm o}}, \;
\phi_{{\rm o}} = P_{{\rm o}} J_A \psi_{{\rm e}}.
$$
Then
$$
\langle [(J - aI)^2 - a^2 I] \psi, \psi \rangle = \langle ( B^2 -
2aB ) \psi_{{\rm o}} , \psi_{{\rm o}} \rangle + \| \phi_{{\rm e}}
\|^2 + \| \phi_{{\rm o}} \|^2 + 2 {\rm Re } \langle ( B - 2a )
\phi_{{\rm o}} , \psi_{{\rm o}} \rangle.
$$
In view of \eqref{step1}, \eqref{lemmagoal} follows from this once
we show that
\begin{equation}\label{step2}
\left| \langle ( B - 2a ) \phi_{{\rm o}} , \psi_{{\rm o}} \rangle
\right| \le \frac12 \left[ \langle ( B^2 - 2aB ) \psi_{{\rm o}} ,
\psi_{{\rm o}} \rangle + \| \phi_{{\rm e}} \|^2 + \| \phi_{{\rm
o}} \|^2 \right ].
\end{equation}
Since
$$
\left| \langle ( B - 2a ) \phi_{{\rm o}} , \psi_{{\rm o}} \rangle
\right| \le \| \phi_{{\rm o}} \| \cdot \| ( B - 2a ) \psi_{{\rm
o}} \| \le \frac12 \| \phi_{{\rm o}} \|^2 + \frac12 \cdot \| ( B -
2a ) \psi_{{\rm o}} \|^2,
$$
\eqref{step2} follows from
\begin{equation}\label{step3}
\frac12 \cdot \| ( B - 2a ) \psi_{{\rm o}} \|^2 \le \langle ( B^2
- 2aB ) \psi_{{\rm o}} , \psi_{{\rm o}} \rangle.
\end{equation}
But \eqref{step3} is a consequence of
$$
0 \le \langle ( B^2 - 4a^2 ) \psi_{{\rm o}} , \psi_{{\rm o}}
\rangle
$$
which in turn follows immediately from the assumption $\psi \in
\mathcal{H}^{(2)}_a$.
\end{proof}

The following theorem describes the positive spectrum of $J$,
which turns out to be discrete in $[0,\infty)$. It is therefore
interesting to study the distribution of eigenvalues
$$
0 < E_1(b) < E_2(b) < \cdots
$$
Note that all eigenvalues must be simple. Let $N(E)$ denote the
number of eigenvalues of $J$ in the interval $(0,E)$.

\begin{theorem}
Let $0 < \alpha \le 1$ and $b > 0$.
\\
{\rm (a)} $\sigma(J) \cap (0,\infty)$ is discrete and can
accumulate only at $\infty$, not at $0$.
\\
{\rm (b)}  The eigenvalue counting function $N$ obeys the estimate
$$
N(E) \le \frac12 \left( \frac{E}{b} \right)^{\frac{1}{\alpha}}.
$$
{\rm (c)} The $n$-th eigenvalue, $E_n(b)$, obeys the lower bound
$$
E_n(b) \ge 2^\alpha b n^\alpha.
$$
{\rm (d)} If $b \ge \sqrt{6}$, then
$$
E_n(b) \le \frac{2^{2 + \alpha}}{(3^{1/\alpha} - 1)^\alpha} b
n^\alpha.
$$
{\rm (e)} If $b < \sqrt{6}$, then
$$
E_n(b) \le C b^{1 - 2\alpha} n^\alpha
$$
for some constant $C > 0$.
\end{theorem}

\begin{proof} (a) and (b) are immediate consequences of
Lemma~\ref{l:keyest}: On $\mathcal{H}^{(2)}_a$, we have $(J-aI)^2
\ge a^2$, and hence the dimension of the range of the spectral
projection of $J$ onto the interval $(0,2a)$ is bounded by
$$
\text{dim} \mathcal{H}^{(1)}_a = \frac12 \left( \frac{2a}{b}
\right)^{1/\alpha}.
$$
(c) follows from (b): For $\varepsilon > 0$, we have
$$
n = N(E_n(b) + \varepsilon) \le \frac12 \left( \frac{E_n(b) +
\varepsilon}{b} \right)^{1/\alpha},
$$
and hence
$$
E_n(b) + \varepsilon \ge 2^\alpha b n^\alpha.
$$
Now let $\varepsilon$ go to zero. This yields the claimed lower
bound for $E_n(b)$.

In order to prove the upper bounds in (d) and (e), we note the following. If, for $n \ge
1$ and $a>0$, we can find test functions $f^{(1)}, \ldots, f^{(n)}$ with disjoint
supports obeying the estimate
\begin{equation}\label{energy}
\|(J - aI)f^{(k)}\|^2 \le a^2 \|f^{(k)}\|^2
\end{equation}
for $1 \le k \le n$, then $E_n(b) \le 2a$.

We first prove the estimate in (d). Our test functions will be of
the form $\delta_{2m-1}$, where $m$ is chosen such that
$$
|b (2m-1)^\alpha - a| < \frac{a}{2}.
$$
Equivalently,
\begin{equation}\label{interval}
\left( \frac{a}{2b} \right)^{1/\alpha} < 2m-1 < \left(
\frac{3a}{2b} \right)^{1/\alpha}.
\end{equation}
For these test functions, we have
\begin{align*}
\|(J-aI)\delta_{2m-1}\|^2 - a^2 \| \delta_{2m-1} \|^2 & \le 2
(2m-1)^{2\alpha}
+ |b(2m-1)^\alpha - a|^2 - a^2\\
& < \frac{9a^2}{2b^2} + \frac{a^2}{4} - a^2.
\end{align*}
If $b \ge \sqrt{6}$, this shows that $\|(J-aI)\delta_{2m-1}\|^2 -
a^2 \| \delta_{2m-1} \|^2$ is negative. Note that the number of
test functions we obtain this way, $n$, is restricted by the
condition
$$
n \le \frac12 \left[ \left( \frac{3a}{2b} \right)^{1/\alpha} -
\left( \frac{a}{2b} \right)^{1/\alpha} \right].
$$
Thus, if we choose
$$
a \ge 2b \left( \frac{2n}{3^{1/\alpha} - 1} \right)^\alpha,
$$
we can find $n$ test functions with disjoint supports obeying
\eqref{energy}. This yields the estimate in (d).

Let us turn to part (e). Again, our test functions will be
supported on odd sites $2m-1$, subject to \eqref{interval}. We
partition the interval in \eqref{interval} into $n$ subintervals
$J_1, \ldots, J_n$ of equal size, roughly given by
$$
\Delta_n = \frac{1}{n} \left[ \left( \frac{3a}{2b} \right)^{1/\alpha} - \left(
\frac{a}{2b} \right)^{1/\alpha} \right].
$$
Since we need this to be at least $2$, we get a first condition for $a$:
\begin{equation}\label{condition0}
\Delta_n \ge 2.
\end{equation}
We will return to this condition below.

The function $f^{(k)}$ will be supported on the odd sites within $J_k$ and on these
sites, it alternates between the values $1$ and $-1$. Our goal is to prove the estimate
\eqref{energy}. By definition of $f^{(k)}$, we have that
$$
\| f^{(k)} \|^2 \approx \frac{\Delta_n}{2}.
$$
Using the fact that $J_k$ lies within the interval \eqref{interval}, we obtain
\begin{align*}
\| (J - aI) f^{(k)} \|^2 & \le \Delta_n \cdot \sup_{2m-1 \in J_k} |b (2m-1)^\alpha - a|^2
+ \mathrm{const} \cdot \sup_{2m-1 \in J_k} (2m)^{2\alpha} + \\
& \qquad + \sum_{2m-1 \in J_k} | (2m+1)^\alpha - (2m-1)^\alpha |^2 \\
& \le \Delta_n \frac{a^2}{4} + \mathrm{const} \cdot \left[ \left( \frac{3a}{2b}
\right)^{1/\alpha} \right]^{2\alpha} + \mathrm{const} \cdot \Delta_n \sup_{2m-1 \in J_k}
(2m-1)^{2(\alpha - 1)}
\end{align*}
with suitable uniform constants. Thus, in order to satisfy \eqref{energy}, we need
\begin{equation}\label{remainder}
\mathrm{const} \cdot \left[ \left( \frac{3a}{2b} \right)^2 + \Delta_n \sup_{2m-1 \in J_k}
(2m-1)^{2(\alpha - 1)} \right] \le \Delta_n \frac{a^2}{4}.
\end{equation}
Let us first prove
\begin{equation}\label{remainder1}
\mathrm{const} \cdot \left( \frac{3a}{2b} \right)^2  \le \Delta_n \frac{a^2}{8}.
\end{equation}
This inequality is equivalent to
\begin{equation}\label{condition1}
\left( \frac{18 \cdot \mathrm{const}}{\left( \frac32 \right)^{1/\alpha} - \left( \frac12
\right)^{1/\alpha} }\right)^\alpha b^{1-2\alpha} n^\alpha \le a.
\end{equation}
Note that the factor in front of $b^{1-2\alpha}$ is uniformly bounded from above by some
$C > 0$ for all $\alpha \in (0,1]$.

Since, by \eqref{interval}, $2m-1 \ge \left( \frac{a}{2b} \right)^{1/\alpha}$ on $J_k$,
we have
$$
\Delta_n \sup_{2m-1 \in J_k} (2m - 1)^{2(\alpha - 1)} \le \Delta_n \left[ \left(
\frac{a}{2b} \right)^{1/\alpha} \right]^{2(\alpha - 1)}
$$
Thus, to complement \eqref{remainder1}, it suffices to show
\begin{equation}\label{remainder2}
\mathrm{const} \cdot \Delta_n \left( \frac{a}{2b} \right)^{\frac{2(\alpha - 1)}{\alpha}}
\le \Delta_n \frac{a^2}{8}
\end{equation}
This inequality is equivalent to
\begin{equation}\label{condition2}
(8 \cdot \mathrm{const})^{\alpha/2} (2b)^{1-\alpha} \le a.
\end{equation}
Since we already have the condition $a \ge C \cdot b^{1-2\alpha}$ from
\eqref{condition1}, \eqref{condition2} is automatically satisfied if $C$ obeys, in
addition to the property above,
$$
C \ge \max_{0 < \alpha \le 1} 8 \cdot \mathrm{const}^{\alpha/2} 2^{1-\alpha}
(\sqrt{6})^\alpha.
$$
(Here we used $b \le \sqrt{6}$.)

Finally, we need to satisfy \eqref{condition0}. Using $a \ge C b^{1-2\alpha} n^\alpha$,
one checks that $\Delta_n \ge 2$ holds as long as
$$
\frac{C \left[ \left( \frac32 \right)^{1/\alpha} - \left( \frac12 \right)^{1/\alpha}
\right]^\alpha}{2^\alpha (\sqrt{6})^{2\alpha}} \ge 1.
$$
(Here we used again that $b \le \sqrt{6}$.) This can again be satisfied uniformly in
$\alpha \in (0,1]$ by choosing $C$ large enough.

Putting everything together, if $C$ is so large that it satisfies the three properties
found above, then for
$$
a \ge C b^{1-2\alpha} n^\alpha,
$$
we can find $n$ test functions $f^{(1)},\ldots,f^{(n)}$ obeying \eqref{energy}.
\end{proof}

%
%
%
%

\section{Open Problems}\label{six}

In this section we list a number of questions and directions for
future research that are suggested by our and previous work.

\begin{enumerate}

\item More general entries: We expect that one can relax the
assumption $\mu_n = r_n = n^\alpha$ and perform a similar
analysis, akin to \cite{jn2}, for more general unbounded
sequences.

\item The $n$-dependence of $E_n(b)$: In fact, a stronger
statement on the behavior of $E_n(b)$ for large $n$ would be
desirable. Is it possible to find a function $F(b)$ such that
$E_n(b) \sim F(b) n^\alpha$ (in the sense that $E_n(b) F(b)^{-1}
n^{-\alpha} \to 1$ as $n \to \infty$)?

\item The $b$-dependence of $E_n(b)$: Some improved estimates on the dependence of
$E_n(b)$ on $b$ would be of interest; especially in the case $\alpha \in (1/2,1)$ for
which the statement in Theorem~5.(e) does not appear to be optimal.

\item Eigenfunction asymptotics for positive energies: We established explicit
asymptotics of the solutions to the associated difference equation for all energies in
$(-\infty,0]$. Is it possible to prove similar results for energies in $(0,\infty)$? We
would expect the formula from Theorem~\ref{negensolas} to hold also for positive
energies.

\item More detailed questions as $\tilde{b}$ crosses zero:
Combining the results of \cite{jn2} and the present paper, we
obtain a rather detailed picture of the spectrum and the spectral
type as $\tilde{b}$ makes a transition through the critical value
$\tilde{b} = 0$ (cf.\ the discussion at the end of Section~1).
This suggests that even more detailed questions could be
addressed. For example, what happens to the spectral function as
$\tilde{b}$ passes through zero? Do we see some concentration at
the eventual eigenvalues?

\end{enumerate}


\begin{thebibliography}{10}

\bibitem{d} J.\ Dombrowski, Absolutely continuous measures for systems of
orthogonal polynomials with unbounded recurrence coefficients,
\textit{Constr.\ Approx.} {\bf 8} (1992), 161--167

\bibitem{dp1} J.\ Dombrowski and S.\ Pedersen, Absolute continuity for
unbounded Jacobi matrices with constant row sums, \textit{J.\
Math.\ Anal.\ Appl.} {\bf 267} (2002), 695--713

\bibitem{dp2} J.\ Dombrowski and S.\ Pedersen, Spectral transition
parameters for a class of Jacobi matrices, \textit{Studia Math.}
{\bf 152} (2002), 217--229

\bibitem{e2} J.\ Edward, Spectra of Jacobi matrices, differential equations on
the circle, and the $su(1, 1)$ Lie Algebra, \textit{SIAM J.\
Math.\ Anal.} {\bf 34} (1993), 824--831

\bibitem{e} S.\ N.\ Elaydi, \textit{An Introduction to Difference Equations},
2nd edition, Undergraduate Texts in Mathematics, Springer-Verlag,
New York (1999)

\bibitem{Eastham} M.\ S.\ P.\ Eastham, \textit{The Asymptotic Solution of
Linear Differential Systems. Applications of the Levinson
Theorem}, London Mathematical Society Monographs. New Series, {\bf
4}, Oxford Science Publications, The Clarendon Press, Oxford
University Press, New York (1989)

\bibitem{gp} D.\ J.\ Gilbert and D.\ B.\ Pearson, On subordinacy and analysis of
the spectrum of one-dimensional Schr\"odinger operators, \textit{J.\ Math.\
Anal.\ Appl.} {\bf 128} (1987), 30--56

\bibitem{jm} J.\ Janas and M.\ Moszy\'nski, Spectral properties of Jacobi
matrices by asymptotic analysis, \textit{J.\ Approx.\ Theory} {\bf
120} (2003), 309--336

\bibitem{jn1} J.\ Janas and S.\ Naboko, Spectral properties of selfadjoint Jacobi
matrices coming from birth and death processes, in \textit{Recent
advances in operator theory and related topics {\rm (}Szeged,
1999{\rm )}}, Oper.\ Theory Adv.\ Appl.\ {\bf 127}, Birkh\"auser,
Basel (2001), 387--397

\bibitem{jn2} J.\ Janas and S.\ Naboko, Spectral analysis of selfadjoint Jacobi matrices
with periodically modulated entries, \textit{J.\ Funct.\ Anal.}
{\bf 191} (2002), 318--342

\bibitem{jn3} J.\ Janas and S.\ Naboko, Criteria for semiboundedness in a class of
unbounded Jacobi operators (Russian), \textit{Algebra i Analiz} {\bf 14} (2002),
158--168; ; translation in \textit{St.\ Petersburg Math.\ J.} \textbf{14} (2003),
479--485

\bibitem{jns} J.\ Janas, S.\ Naboko, and G.\ Stolz, Spectral theory for a class of
periodically perturbed unbounded Jacobi matrices: elementary methods, \textit{J.\
Comput.\ Appl.\ Math.} \textbf{171} (2004), 265--276

\bibitem{jl} S.\ Jitomirskaya and Y.\ Last, Power-law subordinacy and singular spectra.
I. Half-line operators, \textit{Acta Math.} {\bf 183} (1999), 171--189

\bibitem{kp} S.\ Khan and D.\ B.\ Pearson, Subordinacy and spectral theory for
infinite matrices, \textit{Helv.\ Phys.\ Acta} {\bf 65}  (1992),
505--527

\bibitem{lln1} C.\ F.\ Lo, K.\ L.\ Liu, and K.\ M.\ Ng, The multiquantum
intensity-dependent Jaynes-Cummings model with the counterrotating
terms, \textit{Physica A} {\bf 265} (1999), 557--564

\bibitem{lln2} C.\ F.\ Lo, K.\ L.\ Liu, and K.\ M.\ Ng, Exact eigenstates of the
intensity-dependent Jaynes-Cummings model with the
counter-rotating term, \textit{Physica A} {\bf 275} (2000),
463--474

\bibitem{mr} D.\ R.\ Masson and J.\ Repka, Spectral theory of Jacobi matrices in
$\ell\sp 2(\Z)$ and the ${\rm su}(1,1)$ Lie algebra, \textit{SIAM
J.\ Math.\ Anal.} {\bf 22} (1991), 1131--1146

\bibitem{m} M.\ Moszy\'nski, Spectral properties of some Jacobi matrices
with double weights, \textit{J.\ Math.\ Anal.\ Appl.} {\bf 280}
(2003), 400--412

\end{thebibliography}
\end{document}